\documentclass{amuc}

\usepackage{amsfonts}
\usepackage{amsmath}

\begin{document}

%% \title[HEADLINE TITLE]{LONG TITLE \\ WITH TWO LINES}

\title[Non-uniqueness for a nonlocal elliptic equation]{Lack of uniqueness for an elliptic equation with nonlinear and nonlocal drift posed on a torus}

\author[G. Rui]{Giulia Rui}
\address{Department of Mathematics, University of Milano, Italy}
\email{Giulia.Rui@unimi.it}

\author[A. Muntean]{Adrian Muntean}
\address{Department of Mathematics and Computer Science, University of Karlstad, Sweden}
\email{Adrian.Muntean@kau.se}


{35B32}

\begin{abstract}
We study a nonlinear and nonlocal elliptic equation posed on the flat torus.
While constant solutions always exist, we show that uniqueness fails in general.
Using spectral analysis and the Crandall--Rabinowitz bifurcation theorem, we prove
the existence of branches of non-constant periodic solutions bifurcating from
constant states. This result is qualitative and non-constructive. Using a conceptually different argument, we construct explicit multiple solutions for a specific one--dimensional formulation of our target problem.
\end{abstract}

\maketitle

%%%%% Options of setting environments for theorems, definitions, etc.
\theoremstyle{plain} %% for italic environments
\newtheorem{theorem}{Theorem}[section]
\newtheorem{corollary}[theorem]{Corollary}
\newtheorem{lemma}[theorem]{Lemma}

\theoremstyle{definition} %% for nonitalic environments
\newtheorem{definition}[theorem]{Definition}
%%% for unnumbered environments, use:
\newtheorem*{remark}{Remark}

%%%%%%%%%%%%%%%%%%%%%%%%%%%%%%%%%%%%%%%%%%%%%%%%%%%%%%%%%%%%%

\section{Introduction}

In this Note, we investigate  the following elliptic boundary-value problem with nonlinear and nonlocal drift posed on the torus $\mathbb{T}^d$, namely 
\[
\nabla\cdot\bigl(a\nabla u + u\,\Phi*\nabla u\bigr)=0
\qquad \text{in } \mathbb{T}^d,\quad\quad\quad\quad (\star)
\]
where $a>0$ and $\Phi$ is a periodic interaction kernel.
This type of boundary-value problem arises in the context of periodic homogenization of parabolic problems that have nonlocal drifts. Depending on the scaling of the potential $\Phi$, such formulation $(\star)$ appears as a cell problem in the said context\footnote{It is worth noting that in the homogenization context, cell problems are mathematical objects typically used to incorporate microscopic information into effective coefficients.}.

Our motivation to study  ($\star$) stems from exploring a certain class of stochastically interacting particle systems and their mean-field limits. Concretely, in some applications one considers the evolution of $N$ particles whose motion is biased by pairwise interactions through a potential (or kernel) $\Phi$, which, at the hydrodynamic level, leads to a drift velocity given by a convolution field of the form $\Phi * \nabla u$, where $u$ denotes the particle density. In this perspective, the term $u(\Phi*\nabla u)$ arising in $(\star)$ represents a macroscopic signature of possible aggregation and repulsion originally taking place at microscopic scales. 

Somewhat similar nonlocal interaction mechanisms appear across fields, including collective behavior, swarming, granular media, opinion dynamics, and McKean--Vlasov-type models \cite{Carrillo2003Kinetic,Bolley2013Uniform}.
In the specific setting we have in mind, the problem ($\star$) arises in the context of multiscale models of sulphation reactions in marble (or in some other porous material), which are important for cultural heritage conservation \cite{Javergaard2026Hybrid}: The spatial distribution of reactive agents is influenced by interactions such as effective attraction/repulsion induced by the material's microstructure or by eventual crowding effects at the pore level.

The solvability question for the boundary-value problem $(\star)$ can be answered in a straightforward manner. While constant solutions are always present, we wonder about their uniqueness within a fixed-mean class. In this sense, the bifurcation result below
identifies a genuinely nontrivial form of non-uniqueness. This is precisely where this Note contributes: we indicate that the lack of uniqueness of solutions is somewhat inherent in this context. Our main result states that, under mild spectral assumptions on the potential $\Phi$, non-constant periodic solutions bifurcate from constant states while preserving the mean, naturally yielding non-uniqueness of solutions to the target problem. The Crandall--Rabinowitz bifurcation theorem will be the main tool used here to clarify the matter.
For a long-time behaviour and bifurcation analysis of aggregation-diffusion systems on the torus, see also the recent related work 
\cite{CarrilloSalmaniw2025}.

\section{An observation and the statement of our  main result}

\subsection{Observation}
For  convenience, we let $\Omega:=\mathbb{T}^d$ be the flat $d-$dimensional torus. 
We consider the equation
\begin{equation}\label{eq:main}
\nabla\cdot\bigl(a\nabla u + u\,\Phi*\nabla u\bigr)=0
\qquad \text{in } \Omega,
\end{equation}
where $a>0$, $\Phi:\Omega\to\mathbb R $ is an even function $\Phi \in L^2(\Omega)$, and the convolution is taken in the standard sense.

From the modeling perspective, one is often primarily interested in
nonnegative interaction kernels. However, the mathematical mechanism leading
to non-uniqueness does not rely on the sign of $\Phi$, but rather on the
spectral condition imposed below on its Fourier coefficients. We therefore
state the result for real-valued even kernels, which slightly extends the
class suggested by the original modeling motivation.

We note with ease  that
any constant function $u\equiv c$ is a solution of \eqref{eq:main}. Thinking of the uniqueness question, it is interesting to observe that even in the corresponding linear but nonlocal situation, the uniqueness of solutions up to a constant is not guaranteed. To see this fact a bit closer, we first consider a linearized version of $(\star)$, that is, for any given $b\in \mathbb{R}$, we now look at periodic solutions to 
\begin{equation}\label{eq:linear}
a\Delta u + b\Phi*\Delta u = 0
\qquad \text{in } \Omega.
\end{equation}

Throughout the paper, we use the Fourier convention on $\mathbb T^d=\mathbb R^d/\mathbb Z^d$
given by
\[
    \widehat f(k)=\int_{\mathbb T^d} f(y)e^{-2\pi i k\cdot y}\,dy,
    \qquad
    f(y)=\sum_{k\in\mathbb Z^d}\widehat f(k)e^{2\pi i k\cdot y}.
\]
With this normalization, convolution satisfies
\[
    \widehat{\Phi*f}(k)=\widehat\Phi(k)\widehat f(k).
\]

Rewriting \eqref{eq:linear} using Fourier series, we get
\[
u(y)=\sum_{k\in\mathbb Z^d}\widehat u(k)e^{2\pi i k\cdot y},
\qquad
\Phi(y)=\sum_{k\in\mathbb Z^d}\widehat \Phi(k)e^{2\pi i k\cdot y},
\]
which directly leads to
\[
-(2\pi)^2|k|^2\bigl(a+ b\widehat \Phi(k)\bigr)\widehat u(k)=0.
\]
Hence, non-constant solutions exist if and only if there exists some $k\neq0$ such that $\widehat \Phi(k)=-a/b$. In this linear case, it is easy to obtain some conditions on the kernel such that uniqueness can be granted in an {\em a priori} fashion: either via energy estimates, requiring that $\|\Phi\|_{L^1} < a/|b|$, or by imposing  that the Fourier coefficients of the kernel satisfy $\widehat \Phi(k) \neq - a/b$.
If, for instance,  $b>0$, then all kernels with nonnegative Fourier coefficients would guarantee the uniqueness of the solution. In the nonlinear case, however, the situation becomes more complicated, and, as we will see next, for a very wide class of kernels, the uniqueness property can be actively disproved.

\subsection{Statement of the main result}

Let $s>d/2$. We define the spaces
\begin{align}\label{def:spaces}
    X &:= \Bigl\{ v\in  H^{s+2}(\Omega) :
\int_{\Omega} v=0,\ v(y)=v(-y)\Bigr\},\\
Y &:= \Bigl\{ w\in H^{s}(\Omega) :
\int_{\Omega} w=0,\ w(y)=w(-y)\Bigr\}.
\end{align}

We next state the main result of this Note:
\begin{theorem}[Non-uniqueness of periodic solutions to $(\star)$]\label{thm:nonuniqueness}
Assume that there exists an index $k_0\in\mathbb Z^d\setminus\{0\}$ such that
\begin{equation}\tag{H}\label{Kernel_nondegeneracy}
    \widehat \Phi(k_0)\neq0,
    \qquad
    \widehat \Phi(k)\neq \widehat \Phi(k_0)
    \quad\text{for all } k\neq\pm k_0.
\end{equation}

Define $\mathcal F:\mathbb R\times X\to Y$ as
\[
\mathcal F(c,v)
:=
a\Delta v+\nabla\cdot\bigl((c+v)\,\Phi*\nabla(c+v)\bigr),
\]
and the point
\[
    c_0=-\frac{a}{\widehat\Phi(k_0)}.
\]

Then there exist $\varepsilon>0$, a neighborhood $U$ of $(c_0,0)$ in
$\mathbb R\times X$, and  two $\mathcal{C}^1$ curves $\varphi, \psi$ such that
\[
\mathcal F^{-1}(0)\cap U
=
\bigl\{(\varphi(s),\,s v_0+s\psi(s)):\ |s|<\varepsilon\bigr\}
\cup
\bigl\{(c,0): |c-c_0|<\varepsilon\bigr\}.
\]
for a suitable $v_0 \in \mathbb{H}^{s+2}(\Omega)$, $s> d/2$. Specifically, there exists a differentiable curve of nonconstant solutions to \eqref{eq:main}.
\end{theorem}

\begin{remark}
It is worth noticing that the non-uniqueness obtained above occurs within a
fixed-mean class. Indeed, for $t\neq0$ sufficiently small,
\[
    u_t=\varphi(t)+t v_0+t\psi(t)
\]
is a nonconstant solution of \eqref{eq:main}. Since $v_0$ and $\psi(t)$ belong
to $X$, the perturbation $t v_0+t\psi(t)$ has zero mean. Hence $u_t$ and the
constant solution $u\equiv\varphi(t)$ have the same mean, and therefore the
same mass.
\end{remark}
\begin{remark}
If the solutions are interpreted as particle densities, then one should require
the bifurcation to occur from a positive constant state. This amounts to the additional condition
$\widehat\Phi(k_0)<0$. Without it, the bifurcation may occur
from a negative constant state, which is mathematically admissible but not
consistent with a density interpretation.
\end{remark}

Essentially, this tells us that the uniqueness of solutions is, in fact, a property that is not expected to hold true.  
In the next section, we provide a proof of this result.

\section{Proof of Theorem~\ref{thm:nonuniqueness}}
The main tool we use is contained in Lemma 1.1 of \cite{CrandallRabinowitz1973}; see also  a more general version of this result in Theorem 1.7 of \cite{CRANDALL1}. For the reader's convenience and expositional clarity, we recall here this key tool. 
\begin{theorem}[Crandall--Rabinowitz]\label{thm:CR}
Let $X,Y$ be Banach spaces, $V\subset X$ an open neighborhood of $0$, and
$\mathcal F\in C^2(({-1},1)\times V;Y)$ in Fr\'echet sense.
Assume:
\begin{enumerate}
\renewcommand{\labelenumi}{(\roman{enumi})}
\item $\mathcal F(\lambda,0)=0$ for all $\lambda\in(-1,1)$;
\item $\ker D_u\mathcal F(\lambda_0,0)=\mathrm{span}\{x_0\}$ with $x_0\neq0$;
\item $\mathrm{codim}\,\mathrm{Ran}(D_u\mathcal F(\lambda_0,0))=1$;
\item $D_{\lambda u}\mathcal F(\lambda_0,0)[x_0]\notin\mathrm{Ran}(D_u\mathcal F(\lambda_0,0))$.
\end{enumerate}
Then $(\lambda_0,0)$ is a bifurcation point of $\mathcal F(\lambda,u)=0$.
More precisely, if $Z$ is any complement of span $x_0$ in $X$, then there exists a small parameter $\varepsilon>0$ and two continuous maps $\varphi:(-\varepsilon,\varepsilon)\to\mathbb R$,
$\psi:(-\varepsilon,\varepsilon)\to Z$, with $\varphi(0)=\lambda_0, \psi(0)=0$, such that
\[
\mathcal F^{-1}(0)\cap U = \{(\varphi(s),\,s x_0+s\psi(s)):\ |s| <\varepsilon\} \cup \{(\lambda,0):|\lambda|<\varepsilon\}.
\]
If $D_{xx}\mathcal F$ is also continuous, the functions $\varphi,\psi$ are once continuously differentiable.
\end{theorem}

To prove the statement of Theorem~\ref{thm:nonuniqueness}, we are going to show that the conditions of the Crandall--Rabinowitz local bifurcation theorem (i.e. Theorem \ref{thm:CR}) are satisfied.
Informally, this Theorem essentially says that if you have a family of equations $\mathcal F(\lambda,u)=0$ with a trivial branch $u\equiv \lambda$, and the linearized operator loses invertibility at some $\lambda=\lambda_0$ with a $1-$dimensional kernel, then there exists a smooth branch of nontrivial small amplitude  solutions that bifurcates from $(\lambda_0,0)$.

\textbf{Step 1: Regularity and condition (i).}

Recall the Banach spaces $X,Y$ defined  in \eqref{def:spaces}. 
Since $s>d/2$ by hypothesis, $ H^s(\Omega)$ is a Banach algebra (cf. Theorem 4.39 in \cite{AdamsFournier}). 
We recall that our key operator
$\mathcal F: \mathbb R\times X \to Y$ is defined as
\[
    \mathcal{F}(c,v) := a \Delta v + \nabla\cdot( (c+v)\Phi*\nabla(c+v)).
\]
This operator is well-defined, since we benefit from the following standard estimates:
\begin{equation}\label{eq:Hs_estimate}
\begin{aligned}
\|\Delta v\|_{H^s(\Omega)}^2
&\le c_s\,\|v\|_{H^{s+2}(\Omega)}^2,\\
\|\Phi * \Delta v\|_{H^s(\Omega)}
&\le \|\Phi\|_{L^1(\Omega)}\,\|v\|_{H^{s+2}(\Omega)},\\
\|\nabla (v\,\Phi * \nabla v)\|_{H^s(\Omega)}
&= \|\nabla v\,\Phi * \nabla v + v\,\Phi * \Delta v\|_{H^s(\Omega)} \\
&\le C\|\nabla v\|_{H^s(\Omega)}^2 \|\Phi\|_{L^1(\Omega)}
     + C\|v\|_{H^s(\Omega)} \|\Phi\|_{L^1(\Omega)}\|v\|_{H^{s+2}(\Omega)} \\
&\le C\|v\|_{H^{s+1}}^2 \|\Phi\|_{L^1(\Omega)}
     + C\|v\|_{H^s(\Omega)} \|\Phi\|_{L^1(\Omega)}\|v\|_{H^{s+2}(\Omega)} \\
&\le C\|v\|_{H^{s+2}(\Omega)}^2 .
\end{aligned}
\end{equation}

Note that condition $(i)$ is trivially satisfied since constants are always solutions of the original equation.

For the sake of clarity, we recall that a mapping $G:E\to F$ between two Banach spaces $E$ and $F$ is Fr\'echet differentiable (cf. e.g. \cite{Zeidler} Chap. 4.2) at $x\in E$ if there exists a bounded linear operator $DG(x):E\to F$ such that
\[
\lim_{\|h\|_E\to0}\frac{\|G(x+h)-G(x)-DG(x)[h]\|_F}{\|h\|_E}=0.
\]
We are going to show that the mapping $ v \mapsto \mathcal F(c,v)$ can be seen as a polynomial of degree $2$ with  coefficients given by bounded multilinear maps, and that therefore it is infinitely differentiable in the Fr\'echet sense. In fact, we can rewrite $ \mathcal F(c,v)$ as
\[
    \mathcal F(c,v)=a \Delta v + cA v + B(v,v) \qquad \text{for all}\, (c,v) \in \mathbb R \times X
\]
via the operators
\begin{equation*}
\begin{aligned}
A : H^{s+2}(\Omega) \to H^s(\Omega)
    &\qquad& A v &:= \nabla\cdot(\Phi * \nabla v),\\
 B : H^{s+2}(\Omega)\times H^{s+2}(\Omega) \to H^s(\Omega)
    &\qquad& B(v,w) &:= \nabla\cdot\bigl(v(\Phi * \nabla w)\bigr),
\end{aligned}
\end{equation*}
which are respectively linear and bilinear. Both $A$ and $B$ are bounded operators by estimates \eqref{eq:Hs_estimate}.

Explicitly, for any increment $(\delta c,h)\in\mathbb R\times X$, the first Fr\'echet derivative reads
\[
D\mathcal F(c,v)[\delta c,h]
= a\Delta h + \delta c\,A v + c\,A h + B(h,v)+B(v,h),
\]
while the second Fr\'echet  derivative with increments $(\delta c_1,h_1),(\delta c_2,h_2)\in\mathbb R\times X$ is the constant bilinear map
\[
D^2\mathcal F(c,v)[(\delta c_1,h_1),(\delta c_2,h_2)]
=
\delta c_1\,A h_2 + \delta c_2\,A h_1 + B(h_1,h_2)+B(h_2,h_1).
\]

In this case, we take an increment $(\delta c,h)\in\mathbb R\times X$ and compute the expression:
\begin{align*}
\mathcal F(c+\delta c,v+h)-\mathcal F(c,v)
&= a\Delta h
+ \delta c\,\nabla\cdot(\Phi*\nabla v) \\
&\quad + \nabla\cdot\bigl((c+v)\,\Phi*\nabla h\bigr)
+ \nabla\cdot\bigl(h\,\Phi*\nabla v\bigr) \\
&\quad +\delta c\,\nabla\cdot(\Phi*\nabla h) + \nabla\cdot\bigl(h\,\Phi*\nabla h\bigr).
\end{align*}
We define the bounded linear operator $D\mathcal F(c,v):\mathbb R\times H^{s+2}(\Omega)\to H^s(\Omega)$ by
\begin{align*}
D\mathcal F(c,v)[\delta c,h]
&:= a\Delta h
+ \delta c\,\nabla\cdot(\Phi*\nabla v) \\
&\quad + \nabla\cdot\bigl((c+v)\,(\Phi*\nabla h)\bigr)
+ \nabla\cdot\bigl(h\,(\Phi*\nabla v)\bigr).
\end{align*}
Using $s>\frac d2$, the Banach algebra property of $H^{s+1}(\Omega)$ together with the convolution estimate
$\|\Phi*f\|_{H^{s+1}(\Omega)}\le \|\Phi\|_{L^1(\Omega)}\|f\|_{H^{s+1}(\Omega)}$, the terms in the remaining part satisfy
\begin{align*}
\|\delta c\,\nabla\cdot(\Phi*\nabla h)\|_{H^s(\Omega)}
&\le C|\delta c|\,\|h\|_{H^{s+2}(\Omega)},\\
\|\nabla\cdot(h\,(\Phi*\nabla h))\|_{H^s(\Omega)}
&\le C\|h\|_{H^{s+2}(\Omega)}^2.
\end{align*}
Hence, we get
\[
\frac{\|\mathcal F(c+\delta c,v+h)-\mathcal F(c,v)
- D\mathcal F(c,v)[\delta c,h]\|_{H^s(\Omega)}}
{\|(\delta c,h)\|_{\mathbb R\times H^{s+2}(\Omega)}}
\;\xrightarrow[(\delta c,h)\to0]{}\;0,
\]
which proves that $\mathcal F$ is Fr\'echet differentiable.

This also gives as a consequence that the derivative at $v=0$ is
\[
L_c h := D_v\mathcal F(c,0)h
=
a\Delta h + c\,\nabla\cdot(\Phi*\nabla h).
\]

Similarly, the second Fr\'echet derivative with increments $(\delta c_1,h_1),(\delta c_2,h_2)\in\mathbb R\times X$ is a  bilinear map, continuous from $(\mathbb R\times H^{s+2}(\Omega))^2$ into $H^s(\Omega)$:
\begin{align*}
D^2\mathcal F(c,v)[(\delta c_1,h_1),(\delta c_2,h_2)]&=
\delta c_1\,\nabla\cdot(\Phi*\nabla h_2)+ \delta c_2\,\nabla\cdot(\Phi*\nabla h_1) \\
&\quad + \nabla\cdot\bigl(h_1\,\Phi*\nabla h_2\bigr) + \nabla\cdot\bigl(h_2\,\Phi*\nabla h_1\bigr).
\end{align*}
All subsequent derivatives vanish identically, meaning that $\mathcal F$ is infinitely differentiable in the Fr\'echet sense.

\textbf{(ii): The kernel of $L_{c_0}$ has dimension 1.}

Since $L_c$ is translation-invariant on $\mathbb T^d$, it is a Fourier multiplier.
In particular, for each $k\in\mathbb Z^d$, we have
\[
L_c\big(e^{2\pi ik\cdot y}\big)= -(2\pi)^2|k|^2\bigl(a+c\widehat\Phi(k)\bigr)e^{2\pi ik\cdot y} := m(k) e^{2\pi ik\cdot y}.
\]
Hence $\ker L_c$ is spanned by those modes for which $m(k)=0$. Under assumption \eqref{Kernel_nondegeneracy} for $c_0=-\frac{a}{\widehat{\Phi}(k_0)}$  it is nontrivial, and it holds that $m(k)\neq 0$ for all $k\neq \pm k_0$. Since we are restricting ourselves to the even subspace 
\[
    \ker L_{c_0}=\mathrm{span}\{\cos(2\pi k_0\cdot y)\},
\]
the kernel of $L_{c_0}$ is one-dimensional and is spanned by
$v_0(y):=\cos(2\pi k_0\cdot y)$.

\textbf{(iii) The range of $L_{c_0}$ has codimension 1.}
First of all, we characterize the range of $L_{c_0}$. By definition, 
\[
    \mathrm{Ran}(L_{c_0}) := \{ v\in Y: v = L_{c_0}g \,\text{for some}\, g \in X\}.
\]
Since  $L_{c_0}$  is translation invariant, it acts as a Fourier multiplier. Thus $v\in \mathrm{Ran}(L_{c_0})$ if and only if there exists $g\in X$ such that
\begin{equation}\label{eq:rangeLc0}
    \widehat{v}(k) = m(k)\widehat{g}(k), \qquad k \in \mathbb Z^d.
\end{equation}
Since both $v$ and $g$ have zero mean, we must have $\widehat{v}(0)=\widehat{g}(0)=0$.

We intend to prove
\begin{equation}\label{eq:rangechar}
\mathrm{Ran}(L_{c_0})
=
\Bigl\{v\in Y:\ \widehat v(k_0)=\widehat v(-k_0)=0\Bigr\}.
\end{equation}

Take $v\in\mathrm{Ran}(L_{c_0})$. If $k=\pm k_0$, then $m(\pm k_0)=0$. Hence,  equation \eqref{eq:rangeLc0} forces the
compatibility condition
\[
\widehat{v}(\pm k_0)=0.
\]

Conversely, suppose $v\in Y$ satisfies $\widehat v(\pm k_0)=0$.
For $k=0,\pm k_0$, define
\[
    \widehat g(0)=0,\qquad \widehat g(\pm k_0)=0.
\]
For $k\neq 0$ $\pm k_0$, assumption \eqref{Kernel_nondegeneracy} implies $m(k)\neq 0$, which means that any $\widehat{v}(k)$  determines uniquely
\[
\widehat{g}(k)=\frac{\widehat{v}(k)}{m(k)}.
\]
We need to prove that $g\in X$. We first claim that there exists $\eta>0$ such that
\[
    \bigl|a+c_0\widehat\Phi(k)\bigr|\geq \eta
    \qquad \text{for all } k\neq 0,\pm k_0.
\]
Indeed, since $\Phi\in L^2(\mathbb T^d)$, we have
$\widehat\Phi(k)\to0$ as $|k|\to\infty$. Hence
\[
    a+c_0\widehat\Phi(k)\to a>0
    \qquad \text{as } |k|\to\infty.
\]
Thus, the above quantity is bounded away from zero for all sufficiently large
$|k|$. For the remaining modes, there are only finitely many indices, and none of them vanishes by assumption \eqref{Kernel_nondegeneracy}. This proves the claim.

Using such lower bound, together with the fact that for $k\neq0$ one has
\[
    (1+|k|^2)^2\leq C|k|^4,
\]
we obtain
\[
\begin{aligned}
\|g\|_{H^{s+2}}^2
&=
\sum_{k\neq0,\pm k_0}
(1+|k|^2)^{s+2}
\frac{|\widehat v(k)|^2}
{(2\pi)^4} |k|^4|a+c_0\widehat\Phi(k)|^2\\
&\leq
C
\sum_{k\neq0,\pm k_0}
(1+|k|^2)^s|\widehat v(k)|^2
\leq
C\|v\|_{H^s}^2 .
\end{aligned}
\]
Therefore $g\in H^{s+2}$. Since $v$ has zero mean and is even, the function
$g$ also has zero mean and is even. Hence $g\in X$ and $L_{c_0}g=v$.

Therefore, we obtained \eqref{eq:rangechar}.

Since $Y$ consists of even functions, the two conditions
$\widehat v(k_0)=0$ and $\widehat v(-k_0)=0$ are equivalent. Thus \eqref{eq:rangechar} becomes a single linear constraint:
\[
\mathrm{Ran}(L_{c_0}) = \Bigl\{ v\in Y:\ \widehat v (k_0)=0\Bigr\}.
\]
Equivalently, there exists a nonzero continuous linear functional
$\ell:Y\to\mathbb R$ given by $\ell( v):= \widehat v (k_0)$ such that
\[
\mathrm{Ran}(L_{c_0})=\ker \ell.
\]
Hence $\mathrm{Ran}(L_{c_0})$ is a closed hyperplane in $Y$ and
$\mathrm{codim}\,\mathrm{Ran}(L_{c_0})=1$.

\medskip

\textbf{(iv) Transversality condition: $D_c L_c[v_0]\big|_{c=c_0}
\notin \mathrm{Ran}(L_{c_0}).$}

We compute
\[
D_c L_c[v_0]
=
\nabla\cdot(\Phi*\nabla v_0)
=
-(2\pi)^2 |k_0|^2\widehat \Phi(k_0)v_0.
\]
Since $\widehat\Phi(k_0)\neq0$, the Fourier coefficient of
$D_cL_c[v_0]\big|_{c=c_0}$ at $k_0$ is nonzero. On the other hand, by the
range characterization proved above,
\[
    \mathrm{Ran}(L_{c_0})
    =
    \{w\in Y:\widehat w(k_0)=0\},
\]
therefore we have
\[
    D_cL_c[v_0]\big|_{c=c_0}
    \notin \mathrm{Ran}(L_{c_0}).
\]

Now we are in the position to apply Theorem \ref{thm:CR} in this context, which concludes the wanted proof. 

Note that since $\mathcal F(c,v)=a\Delta v+cAv+B(v,v)$ is a polynomial map with bounded multilinear coefficients, its second Fr\'echet derivative is continuous. Therefore, by the last assertion of Theorem~\ref{thm:CR}, the bifurcating curves $\varphi$ and $\psi$ are continuously differentiable.

\begin{remark}
    We expect that Assumption \eqref{Kernel_nondegeneracy} could, in fact, be relaxed. If the dimension of the kernel is higher, then  one may obtain a richer
    local structure of $\mathcal F^{-1}(0)$, possibly containing several curves
    or higher-dimensional sets of solutions -- the non-uniqueness result would still hold, provided that the function $\Phi$ is non constant, and therefore that it has at least one nonzero Fourier coefficient.
    However, this is not proved here;
    such a case would require a more involved argument, for instance via a
    Lyapunov--Schmidt reduction or an equivariant bifurcation analysis.
\end{remark}

\section{A one-dimensional example with explicit construction of a non-trivial solution}

The Crandall-Rabinowitz Theorem does not offer a construction of  an explicit non-constant solution. However, given a specific shape of the kernel, one can sometimes write an explicit non-trivial solution to equation \eqref{eq:main}. This is the aim of this section: we produce such an illustrative example in dimension 1. 
It is worth noting already at this stage that the approach proposed here does not utilize any linearization of \eqref{eq:main}, hence it is conceptually different from the procedure motivated by Theorem \eqref{thm:CR}.

Let $\mathbb T=\mathbb R/\mathbb Z$ be the one-dimensional torus and consider the particular choice of kernel
\[
\Phi(x)=2\cos(2\pi x), \qquad \Phi: \mathbb T \to \mathbb R.
\]
Then the Fourier coefficients of $\Phi$ satisfy
\[
\widehat\Phi(\pm 1)=1,\qquad \widehat\Phi(k)=0 \ \text{for } k\neq \pm 1.
\]
We seek nonconstant solutions of equation
\begin{equation}\label{eq:1dimensionalPb}
        0 = a u_{xx} + (u \Phi*u_x)_x,
\end{equation}
of the form $u(y)=c+v(y)$, where $c\in\mathbb R$ is a constant and $v$ is a periodic, even function with zero mean.

We expand all terms of the equation in the Fourier cosine series. First of all, we note
\begin{equation*}
    v(x) = \widehat{v}(0) + \sum_{n=1}^\infty \widehat{v}(n)\cos(2\pi n x) =  \sum_{n=1}^\infty \widehat{v}(n)\cos(2\pi n x),
\end{equation*}
and then formally we compute:
\begin{eqnarray*}
    v_x(x) &=& - 2\pi \sum n\widehat{v}(n) \sin (2\pi n x)\\
    v_{xx}(x) &=& - (2\pi)^2 \sum n^2\widehat{v}(n) \cos (2\pi n x).
\end{eqnarray*}
Note that the convolution operation corresponds to multiplication in Fourier variables, and, furthermore, $\widehat{\Phi}$ is supported on the modes $\pm 1$. Hence, only these modes contribute.
Consequently, we can write:
\begin{eqnarray*}
    \Phi*v_{xx}(x) &=& - (2\pi)^2\sum n^2\widehat{v}(n)\widehat{\Phi}(n) \cos (2\pi n x) =  - (2\pi)^2\widehat{v}(1)\cos (2\pi x)\\
    \Phi*v_x(x) &=& - 2\pi \widehat{v}(1) \sin (2\pi x)\\
v\Phi*v_x(x) &=& - 2\pi \widehat{v}(1) \sum_{n=1}^\infty \widehat{v}(n)\cos(2\pi n x) \sin (2\pi x)\\
&=& - 2\pi \widehat{v}(1) \sum_{n=1}^\infty \widehat{v}(n) \frac{\sin((n+1)2\pi x) - \sin((n-1)2\pi x)}{2}.
\end{eqnarray*}
Therefore, this yields
\begin{eqnarray*}
        && (v \Phi*v_x)_x (x) = \\
        &=&- 2\pi \widehat{v}(1) \sum_{n=1}^\infty \widehat{v}(n) \pi\left[(n+1)\cos((n+1)2\pi x) - (n-1)\cos((n-1)2\pi x)\right]\\
        &=& - 2\pi^2 \widehat{v}(1)\left[ \sum_{m=2}^\infty \widehat{v}(m-1) m\cos(2\pi m x) - \sum_{m=0}^\infty\widehat{v}(m+1) m\cos(2\pi mx)\right] \\
        &=& - 2\pi^2 \widehat{v}(1) \sum_{m=1}^\infty\left[ \widehat{v}(m-1) - \widehat{v}(m+1) \right] m\cos(2\pi mx),
\end{eqnarray*}
where in the last step we used the zero-mean condition $\widehat v(0)=0$.

For convenience, we adopt the notation $V_n := \widehat{v}(n)$ for the respective coefficients. Finding a solution to \eqref{eq:1dimensionalPb} therefore corresponds to finding all Fourier coefficients via the iteration scheme
\begin{equation*}
 (a+ c \widehat{\Phi}(m))m^2 V_m + \frac{V_1}{2}m\left(V_{m-1} - V_{m+1} \right)  = 0, \qquad m \in \mathbb N^+.
\end{equation*}
The coefficient $V_1=0$ would correspond to a constant solution (because then all the other coefficients would vanish as well). Since we are searching for a non-trivial solution we exclude this case. The recurrence scheme becomes
\begin{equation}\label{eq:recurrence}
\begin{aligned}
V_2 &= 2(a+c), && m=1,\\
V_{m-1}-V_{m+1} &= -\frac{2ma}{V_1}\, V_m,
&& m\in\mathbb N,\; m\ge 2.
\end{aligned}
\end{equation}

Given the structure of the recurrence relation \eqref{eq:recurrence}, we can construct an explicit solution by utilizing the modified Bessel functions of the first kind $\{I_m\}_{m\ge 0}$; see \cite{StegunAbramowitz1972} Chapter 9 for more details, and specifically the item 9.6.10. We recall them here:
\begin{equation}\label{eq:Im_def}
I_m(z)=\sum_{j=0}^\infty \frac{1}{j!\,\Gamma(j+m+1)}\left(\frac{z}{2}\right)^{2j+m},
\qquad z\in\mathbb C,\ m\in\mathbb N\cup 0.
\end{equation}
Interestingly, they satisfy the recurrence \cite{StegunAbramowitz1972} (item 9.6.26):
\begin{equation}\label{eq:Im_recurrence}
I_{m-1}(z)-I_{m+1}(z)=\frac{2m}{z}\,I_m(z),
\qquad m\in \mathbb N, m\ge 1,\ z\neq 0.
\end{equation}
The existence and structure of these functions suggests us to search for a solution of the type
\begin{equation*}
    V_m = C I_m(z), \qquad m \in \mathbb N^+
\end{equation*}
for some constant $C$ to be determined later and some $z\in \mathbb C, z \neq 0$ such that $I_1(z) \neq 0.$ This holds true, for example, for all $z\in \mathbb R^+$. With this choice, by comparing \eqref{eq:recurrence} and \eqref{eq:Im_recurrence} we must have
\begin{equation*}
    - \frac{2ma}{I_1} I_m(z) = C(I_{m-1}- I_{m+1}) = C \frac{2m}{z}I_m(z)
\end{equation*}
which therefore gives
\begin{equation*}
    C = - \frac{az}{I_1(z)}.
\end{equation*}
Moreover, the constant $c$ in \eqref{eq:recurrence} (hence in the structure of the solution we are looking for) can always be chosen such that the recurrence scheme is valid for $m=1$:
\begin{equation*}
    c = \frac{C I_2(z)}{2}-a = - a\left( \frac{z I_2(z)}{2I_1(z)} +1\right).
\end{equation*}

We obtain in this way the explicit family of Fourier coefficients
\begin{equation}\label{eq:explicit_coeffs}
V_m = -\frac{az}{I_1(z)}\,I_m(z)\quad (m\ge 1),
\qquad
V_0=0.
\end{equation}
where $z\in\mathbb R^+$ is a free parameter.

The only thing that remains to be proven is that the computations for the derivatives of $v$ are justified. We therefore need to show that the series 
\begin{equation}\label{eq:explicit_v}
\sum_{m=1}^\infty V_m\cos(2\pi m x) \mbox{ with } 
V_m=-\frac{az}{I_1(z)}\,I_m(z),
\end{equation}
converges uniformly on $\mathbb T$, and therefore defines a smooth, even, periodic function with zero mean, which is a (classical) solution of \eqref{eq:1dimensionalPb}. Moreover, $u$ is obviously non-constant.

From the series definition \eqref{eq:Im_def} with $z\in\mathbb R$ fixed, we have for every $m\ge 0$ the bound
\begin{equation*}
|I_m(z)|
\le \sum_{j=0}^\infty \frac{1}{j!\,(j+m)!}\left(\frac{|z|}{2}\right)^{2j+m}
\le \frac{(|z|/2)^m}{m!}\sum_{j=0}^\infty \frac{(|z|/2)^{2j}}{(j!)^2}.
\end{equation*}
The last series is finite since
\[
    \sum_{j=0}^\infty \frac{(|z|/2)^{2j}}{(j!)^2} \leq \sum_{j=0}^\infty \frac{(|z|/2)^{2j}}{(j!)} = e^{(|z|/2)^2}.
\]
Hence, there exists $C_z>0$ such that
\[
|I_m(z)|\le C_z\frac{(|z|/2)^m}{m!}.
\]
Using \eqref{eq:explicit_coeffs}, this gives
\[
|V_m|\le \frac{a|z|}{|I_1(z)|}\,C_z\frac{(|z|/2)^m}{m!}.
\]
Therefore
\[
\sum_{m=1}^\infty |V_m|
\le \frac{a|z|}{|I_1(z)|}\,C_z\sum_{m=1}^\infty \frac{(|z|/2)^m}{m!}<\infty,
\]
and the cosine series \eqref{eq:explicit_v} converges absolutely and uniformly on $\mathbb T$. Term-wise differentiation is also justified,
since for each $\ell\ge 1$,
\[
\sum_{m=1}^\infty (2\pi m)^\ell |V_m|
<\infty
\]
by the factorial decay of $V_m$. Hence $v\in C^\infty(\mathbb T)$.

\section{Conclusion}

We point out in this Note a situation that naturally leads to \emph{loss of uniqueness} for periodic elliptic problems that are similar to  $\nabla\!\cdot\!\bigl(a\nabla u + u\,\Phi*\nabla u\bigr)=0$ on $\mathbb T^d$.
Essentially, we show that under a mild non-degeneracy condition on a single Fourier mode of the interaction kernel $\Phi$, a local bifurcation argument yields branches of non-constant periodic solutions bifurcating from constant states, thereby guaranteeing non-uniqueness near the bifurcation
point. In addition, for a specific one-dimensional kernel, we constructed an explicit example producing many smooth non-constant solutions.

From the viewpoint of stochastically-interacting particle systems and their corresponding mean-field models, the drift term $u\,(\Phi*\nabla u)$ represents a  self-induced transport generated by pairwise interactions. This Note indicates that, even for smooth periodic interaction kernels, the associated stationary (or cell) problem may admit multiple periodic equilibria. This kind of situation is relevant in several applied settings, including multiscale descriptions of sulphation chemical reactions in porous stones or, seen more broadly, including models used in cultural heritage protection, where the dynamics of reactive agents are impacted by effective attraction or repulsion. A natural next step is to identify selection principles or scalings that pick out physically meaningful equilibria among the possibly multiple steady states.

\begin{acknowledgment}
 We thank Prof. Grigory Panasenko [University of St. Etienne (France) and Vilnius University (Lithuania)] for a fruitful discussion in Karlstad (Sweden) around problem ($\star$). AM thanks the Swedish Research Council (project nr. 2024-05606)  for financial support. G.R. acknowledges membership in GNAMPA (Gruppo Nazionale per l’Analisi Matematica, la Probabilità e le loro Applicazioni) of Istituto Nazionale di Alta Matematica (INdAM), Italy. Last but not least, we are grateful to the referee for the careful reading of the manuscript and for the constructive remarks, which helped us improve the clarity and precision of the paper.
\end{acknowledgment}

\end{document}